\documentclass[a4paper,twoside,10pt]{article}

\usepackage{amsmath}
\usepackage{amssymb}
\usepackage{amsfonts}
\usepackage[english]{babel}
\usepackage{xypic}
\usepackage{theorem}
\usepackage{sectsty}
\usepackage{fullpage}
\usepackage{subfigure}
\usepackage{multirow}
\usepackage{url}

\theorembodyfont{\slshape}
\newtheorem{theo}{Theorem}[section]
\newtheorem{lemm}[theo]{Lemma}
\newtheorem{prop}[theo]{Proposition}
\newtheorem{coro}[theo]{Corollary}
{
 \theorembodyfont{\normalfont}
 \newtheorem{rem}[theo]{Remark}
}

\newenvironment{prf}[1][\!\!]{\noindent\emph{Proof #1 -- }}{\hfill$\square$}

\newcommand\Aa{\mathcal A}
\newcommand\Cc{\mathcal C}
\newcommand\Dd{\mathcal D}
\newcommand\F{\mathbb F}
\newcommand\g{{\mathbf g}}
\newcommand\G{G}
\newcommand\Hh{H}
\newcommand\M{\mathcal M}
\newcommand\Oo{{\mathcal O}}
\newcommand\Pp{\mathbb P}
\newcommand\Q{\mathbb Q}
\newcommand\Ss{\mathcal S}
\newcommand\Z{\mathbb Z}

\newcommand\spl{{\mathrm{split}}}
\newcommand\nspl{{\mathrm{nonsplit}}}

\newcommand\GL{\text{\rm{GL}}}
\newcommand\SL{\text{\rm{SL}}_2}
\newcommand\ST{\text{\rm{ST}}_2}
\newcommand\sla{\text{\rm{sl}}_2}
\newcommand\PSL{\text{\rm{PSL}}_2}
\newcommand\Tr{\text{\rm{Tr\,}}}
\newcommand\Ker{\text{\rm{Ker}}}

\newcommand\z[1]{\Z/{#1}\Z}
\newcommand\gen[1]{{\left\langle #1 \right\rangle}}
\newcommand\vect[2]{\begin{pmatrix}#1\\#2\end{pmatrix}}
\newcommand\svect[2]{\big(\begin{smallmatrix}#1\\#2\end{smallmatrix}\big)}
\newcommand\mat[4]{\begin{pmatrix}#1&#2\\#3&#4\end{pmatrix}}
\newcommand\smat[4]{\big(\begin{smallmatrix}#1&#2\\#3&#4\end{smallmatrix}\big)}

\newcommand\qed{{\hfill\ensuremath{\square}}}

\title{Effective Siegel's Theorem for Modular Curves}
\author{Yuri Bilu\footnote{Universit\'e Bordeaux 1}, Marco Illengo$^\ast$}

1

\begin{document}
\maketitle

\begin{abstract}
We prove that integral points can be effectively determined on all but finitely many modular curves, and on all but one modular curve of prime power level.
\end{abstract}

{\footnotesize
\tableofcontents}

\section{Introduction}
\addtocontents{toc}{\protect\vspace{-0.5\baselineskip}}
Let~$X$ be a curve defined over a number field~$K$ and ${j\in K(X)}$ a $K$-rational function on~$X$. Let also~$R$ be a subring of~$K$. We define the set $X(R,j)$ of $R$-integral points on~$X$ with respect to~$j$ by
$$
X(R,j)=\{P\in X(K): j(P)\in R\}.
$$

The following fundamental theorem was proved by Siegel in 1929.
\begin{theo}[Siegel~\cite{Si29}]
\label{T.siegel}
Assume that either
${\g(X)\geqslant 1}$ or~$j$ has at least~$3$ distinct poles. Then  for any finitely generated subring~$R$ of~$K$ the set ${X(R,j)}$ is finite. 
\end{theo}

Siegel himself considered only the case ${R=\Oo_K}$, but the extension to general~$R$ is relatively straightforward (see \cite{La83,Se97}). Recently Corvaja and Zannier~\cite{CZ02} gave a new beautiful proof of Siegel's theorem, which extends to higher dimensions.  

Theorem~\ref{T.siegel} admits the following converse: if ${\g(X)=0}$ and~$j$ has at most~$2$ distinct poles, then, for some finite extension~$K'$ of~$K$, and some finitely generated subring~$R'$ of~$K'$, the set $X(R',j)$ is infinite. See~\cite{ABP08} for a more precise statement.


For curves of genus at least~$2$, Faltings~\cite{Fa83} improved on the Theorem of Siegel by showing that $X(K)$ is finite if ${\g(X)\geqslant 2}$. 

Both the results of Siegel and Faltings are {\sl non-effective}, that is, neither of them provides any bound for the size of the points in ${X(R,j)}$ computable in terms of~$X$, $j$, $K$ and~$R$.

Let~$X$ be a curve defined over~$\bar\Q$ and ${j\in \bar\Q(X)}$ a non-constant rational function on~$X$. 
We call the couple $(X,j)$ \textsl{Siegelian} if one of the 
conditions of Siegel's Theorem is satisfied, that is: either ${\g(X)\geqslant 1}$ \textsl{or}~$j$ has at least~$3$ distinct poles. Thus, the couple is non-Siegelian if~$X$ is of genus~$0$ \textsl{and}~$j$ has at most two poles. We say that \textsl{Siegel's theorem is effective} for a Siegelian couple $(X,j)$ if for any number field~$K$ such that the couple $(X,j)$ is defined over~$K$, and for any finitely generated subring~$R$ of~$K$, the set $X(R,j)$ (which is finite by the Theorem of Siegel) can be effectively determined in terms of~$X$,~$j$,~$K$ and~$R$.

Starting from pioneering work of A.~Baker, there have been obtained effective versions for some cases of this theorem; see~\cite{Bi95,Bi02} for the history of the subject and further references. For instance, the following is known.

\begin{theo}
\label{T.classic}
Siegel's theorem is effective for $(X,j)$ if
\begin{enumerate}
\addtolength{\itemsep}{-0.5\baselineskip}

\item {\rm(folklore)}
 ${\g(X)=0}$ and~$j$ has at least~$3$ poles, or

\item {\rm(Baker and Coates~\cite{BC71})}
${\g(X)=1}$, or

\item {\rm(Bilu \cite{Bi88}, Dvornicich and Zannier \cite{DZ94})}
${\g(X)\geqslant 1}$ and ${\bar K(X)/\bar K(j)}$ is a Galois extension. 
\end{enumerate}
\end{theo}


Starting from~\cite{Bi95}, Baker's method is applied to obtain effective Siegel's theorem for various classes of modular curves. Let~$\Gamma$ be a congruence subgroup of $\SL(\Z)$ and~$X_\Gamma$ the corresponding modular curve. (See Subsection~\ref{SS.nota} for the definitions.) As usual, we denote~$j$ the modular invariant function. The couple $(X_\Gamma,j)$ is defined over~$\bar\Q$, and one can study the Diophantine properties of this couple. In particular, on can ask the following question:
\begin{quote}
\textsl{assuming the  couple $(X_\Gamma,j)$ Siegelian, is Siegel's theorem effective for this couple?}
\end{quote}

In the sequel, we call a modular curve~$X_\Gamma$ (non-)Siegelian if the couple $(X_\Gamma,j)$ is (non-)Siegelian, where~$j$ is the modular invariant. We shall say that \textsl{Siegel's theorem is effective} for a Siegelian~$X_\Gamma$ if it is effective for the couple $(X_\Gamma,j)$. 

In \cite{Bi95,Bi02} Siegel's theorem was shown to be effective for several classes of modular curves, like the curves $X(n)$, $X_1(n)$ and $X_0(n)$ (provided they are Siegelian). For $X(n)$ effective Siegel's theorem was already established by Kubert and Lang \cite[Section~8.1]{KL81} (they do not make any mention  of effectiveness, but it is implicit in their work).   The results of \cite{Bi95,Bi02} are based on the ``three cusps criterion'', see Section~\ref{S.cusps}.

In the present article we show that Siegel's theorem is effective  for all but finitely many~$X_\Gamma$, and for all but one~$X_\Gamma$ of prime power level. Our principal results are the following two theorems.

\begin{theo}
\label{T.power}
Let~$\Gamma$ be a subgroup of $\SL(\Z)$ of prime power level, distinct from~$25$. Then either~$X_\Gamma$ is non-Siegelian, or Siegel's theorem is effective for~$X_\Gamma$.
\end{theo}

At level~$25$ there is a subgroup~$\Gamma$, defined in Proposition~\ref{P.power_semis}, for which the curve $X_\Gamma$ is of genus~$2$ and for which our argument does not work.

\begin{theo}\label{T.mixed}
Let~$\Gamma$ be a subgroup of level not dividing the number  ${2^{20}\!\cdot\!3^7\!\cdot\!5^3\!\cdot\!7^2\!\cdot\!11\!\cdot\!13}$. Then Siegel's theorem is effective for~$X_\Gamma$.
\end{theo}

The assumption on the level in Theorem~\ref{T.mixed} can certainly be relaxed, but at the moment, the methods of the present article do not allow treatment of certain Siegelian modular curves of small mixed level. Consider, for instance, two congruence subgroups~$\Gamma_5$ and~$\Gamma_7$ of levels~$5$ and~$7$, whose projections to ${\PSL(\F_5)}$ and ${\PSL(\F_7)}$  (see Table~\ref{Tab.p8}) are isomorphic to the fourth alternating group~$\Aa_4$ and to the fourth symmetric group~$\Ss_4$, respectively; their intersection ${\Gamma_5\cap\Gamma_7}$ is a congruence subgroup~$\Gamma$ of the level $35$ such that~$X_\Gamma$ has genus~$2$. This~$X_\Gamma$ is non-Siegelian, but eludes our methods.


\subsection{Notation and Conventions}
\label{SS.nota}
We denote by~$\Cc_n$ the $n$-th cyclic group, and by~$\Dd_n$ the $n$-th dihedral group (so that~$\Cc_n$ is the index~$2$ subgroup of~$\Dd_n$). Further, we denote by~$\Ss_n$ and~$\Aa_n$ the $n$-th symmetric and alternating groups, respectively.

The letter~$\Gamma$ is reserved to congruence subgroups of $\SL(\Z)$, that is, subgroups containing $\Gamma(n)$ for some~$n$. We shall say in this case that~$\Gamma$ is of \textsl{level dividing~$n$}. The smallest~$n$ with this property will be called the \textsl{exact level} of~$\Gamma$. 

Fix a positive integer~$n$. Then to any  congruence subgroup~$\Gamma$ of level dividing~$n$ we associate  a subgroup~$G$ of $\SL(\Z/n\Z)$ and a subgroup~$\bar G$ of $\PSL(\Z/n\Z)$, as the images of~$\Gamma$ under the natural maps ${\SL(\Z)\to\SL(\Z/n\Z)\to\PSL(\Z/n\Z)}$.  Conversely,~$\Gamma$ is uniquely determined by~$G$ (and~$n$). When~$n$ is~$2$,~$4$ or $p^k$ with an odd prime~$p$, then~$\Gamma$ is uniquely determined by~$\bar G$, under the additional assumption that ${\Gamma \ni -I}$. 
In the sequel, when~$n$ is fixed, we shall freely interchange between~$\Gamma$ and~$G$, when it causes no confusion. Also, when the additional assumptions indicated above are satisfied, we shall interchange between~$\Gamma$,~$G$ and~$\bar G$. 

We use the common notation $\nu_2(\Gamma)$, $\nu_3(\Gamma)$, $\nu_\infty(\Gamma)$ and $\mu(\Gamma)$ for, respectively, the number of the $2$-elliptic points of~$\Gamma$, the number of its $3$-elliptic points, the number of its cusps, and the index ${[\PSL(\Z):\bar\Gamma]}$, where~$\bar\Gamma$ is the image of~$\Gamma$ in 
$\PSL(\Z)$. 

The modular curve $X_\Gamma$ is, by definition, the quotient $\Gamma\backslash \bar{\mathcal H}$ (where~$\bar {\mathcal H}$ is the extended Poincar\'e upper half-plane), with properly defined topology and analytic structure. The modular invariant~$j$ defined a non-constant rational function on~$X_\Gamma$, whose poles are exactly the cusps.  While defined analytically, the curve~$X_\Gamma$, or, more precisely, the couple $(X_\Gamma,j)$  has a model over~$\bar\Q$ (even over $\Q(\zeta_n)$, where~$n$ is the level of~$\Gamma$).  See any standard reference like~\cite{La76,Sh71} for all the missing details.

\subsection{Plan of the Article}
In Section~\ref{S.cusps} we state our main tool, the ``three cusps criterion'', and obtain some auxiliary results on the cusps to be used throughout the article. In Section~\ref{S.prime} we study curves of the prime level, and show that for them Siegel's theorem is effective whenever they are Siegelian; we also classify the non-Siegelian curves of prime level. In Section~\ref{S.power} we do the same for the curves of prime power level (with the aforementioned  exception at level~$25$). In Section~\ref{S.mixed} we consider mixed levels.

\paragraph{Acknowledgments}
We thank Roberto Maria Avanzi and Michael Stoll for useful discussions, and  for performing,  at our request,  electronic calculations that helped our intuition. We also thank Pierre Parent for useful discussions.

\section{The ``Three Cusps Criterion''}
\addtocontents{toc}{\protect\vspace{-0.5\baselineskip}}
\label{S.cusps}
The following theorem (see \cite{Bi95}) plays a capital role in the present article.
\begin{theo}
\label{T.B_3cusps}
Let~$\Gamma$ be a congruence subgroup of $\SL(\Z)$. Then Siegel's theorem is effective for~$X_\Gamma$ if 
the group $\Gamma$  has at least $3$ cusps.
\end{theo}

We shall also use the following refinement of Theorem~\ref{T.B_3cusps}, see \cite[Proposition~12]{Bi02}.
\begin{theo}
\label{T.B_subgroup}
Let $\Gamma$ have a congruence subgroup $\Gamma'$, which contains all 
elliptic elements of $\Gamma$ and has at least $3$ cusps.
Then Siegel's theorem is effective for~$X_\Gamma$.
\end{theo}

Applying Theorems~\ref{T.B_3cusps} and~\ref{T.B_subgroup} requires  computing (or estimating) the number of cusps ${\nu_\infty(\Gamma)}$ of a congruence subgroup~$\Gamma$.  For this purpose we shall use the following simple lemma.  It is certainly known, but we could not find a proof in the literature.  

For any natural~$n$ we denote by~$\M_n$ the set of elements of exact order~$n$ in ${\z{n}\times\z{n}}$.
Obviously,
\begin{equation}
\label{E.cardM}
|\M_n|=n^2\prod_{p|n}\left(1-p^{-2}\right),
\end{equation}
the product being taken over all primes~$p$ dividing~$n$.

\begin{lemm}\label{L.cusps}
Let~$\Gamma$ be a congruence subgroup of level dividing~$n$ and containing~${-I}$, and let~$G$ be the projection of~$\Gamma$ modulo~$n$. Then the number ${\nu_\infty(\Gamma)}$ is equal to the number of the orbits of the natural (left) $G$-action on~$\M_n$. In symbols, we have ${\nu_\infty(\Gamma)=|G\backslash\M_n|}$.
\end{lemm}
\begin{prf}
The number ${\nu_\infty(\Gamma)}$ equals the number of $\Gamma$-orbits of ${\Pp_1(\Q)=\Q\cup\{\infty\}}$ and, since~$\Gamma$ contains~${-I}$, is also the number of $\Gamma$-orbits in the set~$\M$ of coprime couples ${(a,b)\in\Z\times\Z}$. It will suffice to prove that~$\M_n$ corresponds to the set of ${\Gamma(n)}$-orbits of~$\M$, where ${\Gamma(n)}$ is the principal congruence subgroup of level~$n$, i.e. the kernel of the reduction map ${\SL(\Z)\twoheadrightarrow\SL(\z{n})}$.

First, let ${\svect{a}{b}\in\M}$ be any representative of ${\svect{1}{0}\in\M_n}$, that is, let~$a$ and~$b$ be any two coprime integers with ${a\equiv1\pmod n}$ and ${b\equiv0\pmod n}$. As~$a$ and~$b$ are coprime, there exist integers~$x$ and~$y$ such that ${ax+by=1}$. Note that ${x\equiv1\pmod n}$. Then the matrix ${M=\smat{x+by}{y-ay}{-b}{a}}$ lies in ${\Gamma(n)}$ and maps $\svect{a}{b}$ to $\svect{1}{0}$.
This shows that the ${\Gamma(n)}$-orbit of $\svect{1}{0}$ is the class of all representative of $\svect{1}{0}\in\M_n$.
We conclude by the transitivity of ${\SL(\Z)}$ over~$\M$ and by the normality of $\Gamma(n)$ in ${\SL(\Z)}$.
\end{prf}

\begin{coro}\label{C.vinfty}
Let~$\Gamma$ and~$G$ be as in the proposition. Assume that~$\Gamma$ has at most~$2$ cusps. Then ${|G|\geqslant|\M_n|/2}$.\hfill\qed
\end{coro}

\section{The Prime Levels}
\addtocontents{toc}{\protect\vspace{-0.5\baselineskip}}
\label{S.prime}
In this section we classify the non-Siegelian modular curves of prime level, and prove effective Siegel's theorem for Siegelian curves of prime level. 

\begin{theo}
\label{T.prime}
\begin{enumerate}\addtolength{\itemsep}{-0.5\baselineskip}

\item
All the~$\Gamma$ (up to conjugacy) of exact prime level, for which $X_\Gamma$ is non-Siegelian, are listed in Tables~\ref{Tab.p8} and~\ref{Tab.pST} on page~\pageref{Tab.p8}.

\item
Let~$\Gamma$ be a congruence subgroup of prime level such that $X_\Gamma$ is Siegelian. Then Siegel's theorem is effective for~$X_\Gamma$. 
\end{enumerate}
\end{theo}

\subsection{Lemmas}
Here we collect basic properties of the special linear group ${\SL(\F_p)}$.
The following property is well-known but we sketch a proof for the sake of completeness. 

\begin{prop}\label{P.order}
The order of an element of ${\SL(\F_p)}$ is either ${2p}$ or at most ${p+1}$. When ${p\neq 2}$, the order of an element of ${\PSL(\F_p)}$ is either~$p$ or at most ${(p+1)/2}$.
\end{prop}
\begin{prf}
A matrix from $\SL(\F_p)$ is either similar over~$\F_p$ to $\smat\lambda10\lambda$ with ${\lambda=\pm1}$ or similar over~$\F_{p^2}$ to $\smat{\alpha}{0}{0}{\alpha^{-1}}$ with ${\alpha\in \F_{p^2}}$. In the first case the order divides ${2p}$. In the second case either ${\alpha\in \F_p}$, in which case the order divides ${p-1}$, or~$\alpha$ is in the kernel of the norm map ${\F_{p^2}\to\F_p}$, in which case the order divides ${p+1}$.
\end{prf}

\medskip

We shall systematically use the classification of semi-simple subgroups of ${\PSL(\F_p)}$. Actually, a classification for ${\mathrm{PGL}_2(\F_p)}$ is available, see \cite[Proposition~16]{Se72}.

\begin{prop}\label{P.serre.np}
Let~$\bar G$ be a proper subgroup of ${\mathrm{PGL}_2(\F_p)}$ of order not divisible by~$p$. Then~$\bar G$ is isomorphic to one of the following groups:
\begin{itemize}\addtolength{\itemsep}{-0.5\baselineskip}

\item $\Cc_n$, the $n$-th cyclic group;
\item $\Dd_n$, the  $n$-th dihedral group;
\item $\Aa_4$, the fourth alternating group;
\item $\Ss_4$, the fourth symmetric group;
\item $\Aa_5$, the fifth alternating group (this only happens when ${p\equiv\pm1\pmod5}$).
\end{itemize}
\end{prop}

In the unipotent case, one has the following, see~\cite[Proposition~15]{Se72}.
\begin{prop}\label{P.serre.p}
Let~$G$ be a subgroup of ${\mathrm{GL}_2(\F_p)}$ of order divisible by~$p$. Then~$G$ either contains ${\SL(\F_p)}$ or is contained in a Borel subgroup of ${\textrm{GL}_2(\F_p)}$.
\end{prop}
(A \textsl{Borel subgroup} of ${\mathrm{GL}_2(\F_p)}$ is a subgroup conjugate to the subgroup ${\mathrm{GT}_2(\F_p)}$ of the upper-triangular matrices.)

\begin{prop}\label{P.st}
Let~$G$ be a subgroup of  the special triangular group ${\ST(\F_p)}$ with ${\nu_\infty(G)\leqslant 2}$.
Then ${G=\ST(\F_p)}$.
\end{prop}
\begin{prf}
If~$G$ were a proper subgroup of ${\ST(\F_p)}$, then its cardinality would be at most half the cardinality of ${\ST(\F_p)}$, that is, 
${|G|\leqslant (p^2-p)/2}$.
On the other hand, 
${|G|\geqslant (p^2-1)/2}$
by Corollary~\ref{C.vinfty}, a contradiction. 
\end{prf}

\begin{theo}\label{T.p.list}
Let~$\Gamma$ be a congruence subgroup of exact level~$p$, with at most~$2$ cusps.
\begin{itemize}\addtolength{\itemsep}{-0.5\baselineskip}

\item If~$p$ does not divide the cardinality of~$\bar G$ then we are in one of the following eight  cases.
\begin{equation}
\label{E.eight}
\begin{tabular}{ll@{\ }l}
 ${p=2}$  & and  ${\bar G\cong \Cc_3}$;\\
 ${p=3}$  & and  ${\bar G\cong \Cc_2}$  & or  $\Dd_2$;\\
 ${p=5}$  & and  ${\bar G\cong \Dd_3}$ & or  $\Aa_4$;\\
 ${p=7}$  & and  ${\bar G\cong \Aa_4}$ & or   $\Ss_4$;\\
 ${p=11}$ & and  ${\bar G\cong \Aa_5}$.
\end{tabular}
\end{equation}
\item If~$p$ divides the cardinality of~$\bar G$ then~$G$ is conjugate to ${\ST(\F_p)}$ and ${\nu_\infty(\Gamma)=2}$.
\end{itemize}
\end{theo}
\begin{prf}
If~${p=2}$ we conclude by inspection. Now assume that ${p\geqslant 3}$. When~$|\bar{G}|$ is not divisible by~$p$, Propositions~\ref{P.order} and~\ref{P.serre.np} imply the upper bound 
${|\bar G|\leqslant \max\{p+1,60\}}$,
and~$60$ can be replaced by~$24$ if ${p\not \equiv \pm1\pmod 5}$. On the other hand, Corollary~\ref{C.vinfty} implies the lower bound 
${|\bar G|\geqslant (p^2-1)/4}$.
It follows that ${p\leqslant 11}$, and we again conclude by inspection. See~\cite[Theorem~6.1.6]{Il08} for more details.

When~$p$ divides~$|\bar{G}|$, Proposition~\ref{P.serre.p} implies that either ${G=\SL(\F_p)}$ or~$G$ is conjugate to a subgroup of ${\ST(\F_p)}$. In the first case~$\Gamma$ is ${\SL(\Z)}$, against our assumption on its level; in the second case we conclude by Proposition~\ref{P.st}.
\end{prf}

\medskip

The invariants of the modular curves corresponding to the eight cases~\eqref{E.eight} are given in Table~\ref{Tab.p8}.
We see that all the corresponding curves are non-Siegelian. We may also remark that in the first five cases (with ${p\leqslant5}$) the group $\bar G$ is uniquely defined up to conjugacy, and that in each of the last three cases (with ${p\geqslant7}$) the group~$\bar G$ belongs to one of two distinct conjugacy classes, so, up to modular equivalence, Table~\ref{Tab.p8} defines~$11$ modular curves. 

Remark that in all the above cases we have $\nu_\infty(G)|G|=|\M_p|$.

\medskip


When~$p$ divides~${|G|}$ and ${\mu_\infty(G)\leqslant2}$, by Theorem~\ref{T.p.list} the either the group~$G$ is  ${\SL(\F_p)}$, in which case we obtain the non-Siegelian curve $X(1)$, or~$G$ is conjugate to ${\ST(\F_p)}$. In this this latter case, up to conjugacy, ${\Gamma=\Gamma_0(p)}$. The effectivity problem for the modular curves $X_0(n)$ is completely solved in \cite[Theorem~10]{Bi02}:

\begin{theo}
Given an integer ${n>1}$, either Siegel's theorem is effective for $X_0(n)$ or the couple $X_0(n)$ is non-Siegelian, which is the case if and only~$n$ is in the set ${\{2,3,5,7,13\}}$.\qed
\end{theo}

The invariants of the corresponding modular curves are given in Table~\ref{Tab.pST}.

This completes the proof of Theorem~\ref{T.prime}. 

\begin{table}[hbt!]
 \caption{Non-Siegelian modular curves of exact prime level}
 \label{Tab.p}
 \subtable[The semi-simple case]{
 \label{Tab.p8}
  \begin{tabular}{cccccccl}
 $p$ & $\bar G$ & $\mu$ & $\nu_\infty$ & $\nu_2$ & $\nu_3$ & $\g$ & remark\\
    \hline\\[-2.5ex]
  2 & $\Cc_3$ &    2 & 1 & 0 & 2 & 0 & \\
  3 & $\Cc_2$ &    6 & 2 & 2 & 0 & 0 &  $X_\spl(3)$\\
  3 & $\Dd_2$ &    3 & 1 & 3 & 0 & 0 & $X_\nspl(3)$\\
  5 & $\Dd_3$ &   10 & 2 & 2 & 1 & 0 & $X_\nspl(5)$\\
  5 & $\Aa_4$ &    5 & 1 & 1 & 2 & 0 & \\
  7 & $\Aa_4$ &   14 & 2 & 2 & 2 & 0 & 2 groups\\
  7 & $\Ss_4$ &    7 & 1 & 3 & 1 & 0 & 2 groups\\
 11 & $\Aa_5$ &   11 & 1 & 3 & 2 & 0 & 2 groups
  \end{tabular}
 }
\hfill
 \subtable[The unipotent case]{
 \label{Tab.pST}
  \begin{tabular}{cccrrrc}
   $p$ & $\Gamma$  & $\mu$ & $\nu_\infty$ & $\nu_2$ & $\nu_3$ & $\g$\\
    \hline\\[-2.5ex]
   2 & $\Gamma_0(2)$   &  3 & 2 & 1 & 0 & 0 \\
   3 & $\Gamma_0(3)$   &  4 & 2 & 0 & 1 & 0 \\
   5 & $\Gamma_0(5)$   &  6 & 2 & 2 & 0 & 0 \\
   7 & $\Gamma_0(7)$   &  8 & 2 & 0 & 2 & 0 \\
  13 & $\Gamma_0(13)$  & 14 & 2 & 2 & 2 & 0 
  \end{tabular}
 }
\end{table}

\section{The Prime Power Levels}
\addtocontents{toc}{\protect\vspace{-0.5\baselineskip}}
\label{S.power}

\subsection{Introduction}
In this section we study groups of prime power level. Our ultimate goal is Theorem~\ref{T.power}. 
As in the prime case, our main tool will be ``three cusps criterion'', in the refined form of Theorem~\ref{T.B_subgroup}

We obtain a complete classification, up to conjugacy, of the groups~$\Gamma$, containing~${-I}$, that do not satisfy the hypothesis of Theorem~\ref{T.B_subgroup}.  Notice that the hypothesis of Theorem~\ref{T.B_subgroup} automatically fails for~$\Gamma$ if~$X_\Gamma$ is non-Siegelian. Thus, as a by-product, we classify non-Siegelian modular curves of prime power level. Up to modular equivalence, there are 34 such curves:
\begin{itemize}
\addtolength{\itemsep}{-0.5\baselineskip}
\item 
the curve $X(1)$ of level~$1$;
\item
16 curves of prime level, listed in Table~\ref{Tab.p};

\item
17 curves of exact level $p^e$ with $e>1$, listed in Table~\ref{Tab.power}.

\end{itemize} 

Besides them, there are three more modular curves of prime power level, for which the hypothesis of Theorem~\ref{T.B_subgroup} is not satisfied. Two of them, one of level~$27$ and the other of level~$32$, are defined in Propositions~\ref{P.power_unip} and~\ref{P.2unipotent} and have genus~$1$; for them Siegel's theorem is effective due to Theorem~\ref{T.classic}. 
The third one, the already mentioned curve of level~$25$ and genus~$2$, occurs in Proposition~\ref{P.power_semis}, and this is the only curve of prime power level for which our argument fails. 


\subsection{The ``Exponential'' Map}

Let~$p$ be a prime number and~$r$,~$s$ positive integers. We denote by $\mathrm M_{2}(R)$ the ring  of ${2\times2}$ matrices over a ring~$R$, and by $\sla(R)$ the additive group of traceless ${2\times2}$ matrices. We define the ``exponential'' map 
$$
\exp=\exp_{r,s}: \mathrm M_{2}(\Z/p^r\Z) \to \GL_2(\Z/p^{s+r}\Z)
$$
by ${\exp(A)=I+p^s\widetilde A}$, where ${\widetilde A\in \mathrm M_{2}(\Z/p^{r+s}\Z)}$ is a lifting of~$A$; clearly, $\exp(A)$ does not depend on the choice of the lifting. Slightly abusing notation, we shall often write ${I+p^s A}$ instead of $\exp(A)$.  

\begin{prop}
Assume that ${r\leqslant s}$. Then  ${\exp\bigl(\sla(\Z/p^r\Z)\bigr)\subset \SL(\Z/p^{r+s}\Z)}$ and we have the short exact sequence
 \begin{equation}\label{E.sequence}
  \xymatrix@1
   {
    \sla(\Z/p^r\Z)\ \ar@{^(->}[r]^-{\exp_{r,s}}&\ \SL(\z{p^{s+r}})\ \ar@{->>}[r]^{\pi_s}&\ \SL(\z{p^s}),
   }
 \end{equation}
 where~$\pi_s$ is the reduction modulo~$p^s$.
\end{prop}

The proof is immediate because 
${\det \bigl(I+p^s\widetilde A\bigr)= 1+ p^s\Tr \widetilde A+ p^{2s}\det \widetilde A}$.

\subsection{Reductions}
Let~$p$ be a prime, ${q=p^e}$ be a power of~$p$, and~$\Gamma$ be a congruence subgroup of exact level~$q$.
For a positive integer~$s$ we consider the \textsl{reduction map} modulo~$p^s$ from ${\SL(\Z)}$ to ${\SL(\z{p^s})}$; the image of~$\Gamma$ is a subgroup~$G_s$ of ${\SL(\z{p^s})}$, whose preimage is a congruence subgroup ${\Gamma_s=\Gamma\cdot\Gamma(p^s)}$ of level dividing~$p^s$.
Then we have a chain of surjective maps
\begin{equation}
\label{E.gamma_g}
 \Gamma\twoheadrightarrow\cdots\twoheadrightarrow G_{e+1}\twoheadrightarrow G_e\twoheadrightarrow G_{e-1}\twoheadrightarrow\cdots\twoheadrightarrow G_2\twoheadrightarrow G_1,
\end{equation}
and a corresponding nested chain of congruence subgroups
\begin{equation*}
 \Gamma=\cdots=\Gamma_{e+1}=\Gamma_e\subsetneq\Gamma_{e-1}\subseteq\cdots\subseteq\Gamma_2\subseteq\Gamma_1.
\end{equation*}
Note that if~$\Gamma$ satisfies the conditions
\begin{equation}
\label{E.-I+lesstwo}
\Gamma\ni -I,\qquad
\nu_\infty(\Gamma)\leqslant 2. 
\end{equation}
then so does~$\Gamma_s$ for every~$s$; in particular, the congruence subgroup~$\Gamma_1$ of level dividing~$p$ belongs to the finite set of groups that we have determined in the previous section.

\begin{rem}
\label{R.es}
One might notice that, while we assume the group ${\Gamma=\Gamma_e}$ to have the exact level~$p^e$, for the ${s\ne e}$ the group~$\Gamma_s$ is not obliged to have the exact level~$p^s$ (actually, it never does for ${s>e}$ and sometimes even for ${s<e}$); \textsl{a priori}, we only know that its level divides~$p^s$. 
\end{rem}

\medskip
For a positive integer~$s$ put ${K_{s+1}=\Ker(\pi_s\vert_{G_{s+1}})}$, where the groups~$G_i$ are as in~\eqref{E.gamma_g} and~$\pi_s$ is the reduction modulo~$p^s$. 
Taking ${r=1}$ in~\eqref{E.sequence} we have a short exact sequence
\begin{equation*}
 \xymatrix@1{
  \sla(\F_p)\ \ar@{^(->}[r]^-{\exp_{1,s}}&\ \SL(\z{p^{s+1}})\ \ar@{->>}[r]^{\pi_s}&\ \SL(\z{p^s}),
 } 
\end{equation*}
and by restriction to the subgroup~$G_{s+1}$ of ${\SL(\z{p^{s+1}})}$ we obtain a short exact sequence
\begin{equation}
 \xymatrix@1{
  V_s\ \ar@{^(->}[r]^-{\exp_{1,s}}&\ G_{s+1}\ \ar@{->>}[r]^{\pi_s}&\ G_s,
 }
\end{equation}
where 
$$
V_s=\exp_{1,s}^{-1}(K_{s+1})
$$ 
is a subspace of ${\sla(\F_p)}$. Thus, the chain of maps~\eqref{E.gamma_g} determines a sequence of subspaces ${V_1,V_2, \ldots}$ of ${\sla(\F_p)}$. 

The group $\SL(\F_p)$ acts by conjugation on $\sla(\F_p)$, which defines a natural action on $\sla(\F_p)$ of any subgroup of $\SL(\F_p)$, in particular of~$G_1$. The following is immediate. 

\begin{prop}
\label{P.Vinvariant}
The spaces $V_s$ are invariant under the  natural action of~$G_1$ on $\sla(\F_p)$ defined above. \qed
\end{prop}

It is crucial that the sequence $(V_i)$ is (non-strictly) increasing, with one little exception. 

\begin{prop}\label{P.Vcontained}
 If ${p^s\neq2}$ then ${V_s\subset V_{s+1}}$. If ${p=2}$ then ${V_1\subset V_2+\gen I}$. 
\end{prop}

\begin{prf}
Let~$M$ be an element of $V_s$, so that~$G_{s+1}$ contains the element ${I+p^sM}$. By surjectivity of the projection ${\pi_{s+1}\colon G_{s+2}\to G_{s+1}}$, there exists a matrix~$N$ with entries in $\z{p^2}$ such that ${I+p^sN\in G_{s+2}}$ projects to  ${I+p^sM}$; obviously, ${N\equiv M\pmod p}$. In~$G_{s+2}$ the~$p$-th power of ${I+p^sN}$ is
\begin{equation*}
 (I+p^sN)^p=I+p^{s+1}N+\vect{p}{2}p^{2s}N^2=I+p^{s+1}\left(N+\vect{p}{2}p^{s-1}N^2\right),
\end{equation*}
implying that ${M+\svect{p}{2}p^{s-1}M^2}$ lies in~$V_{s+1}$. If ${p\neq2}$ or ${s>1}$, then~$p$ divides ${\svect{p}{2}p^{s-1}}$ and therefore ${M\in V_{s+1}}$.
 If ${p=2}$ and ${s=1}$ then ${M+M^2}$ lies in~$V_2$.
Since ${\Tr M=0}$
we have ${M^2=-I\det M}$, whence ${M\in V_2+\gen I}$.
\end{prf}

	\begin{rem}
	\label{R.v_1}
	It is worth mentioning that, when ${p=2}$ and ${-I\in\Gamma}$, we have ${I\in V_1}$, because ${I+2I=-I}$ belongs to~$G_2$. 
	\end{rem}

\begin{coro}\label{C.sla}
Let~$\Gamma$ be a congruence subgroup of the exact level~$p^e$. If ${V_s=\sla(\F_p)}$ for some $s$, then ${e\leqslant s}$.
\end{coro}

\begin{prf}
For ${e>1}$ the hypothesis implies ${V_{e-1}\neq\sla(\F_p)}$. Then it suffices to show that ${V_s=\sla(\F_p)}$ implies ${V_{s+1}=\sla(\F_p)}$. This follows from Proposition~\ref{P.Vcontained} if ${p^s>2}$, and it is verified by inspection for ${p^s=2}$.
\end{prf}

\begin{rem}
 This corollary implies that the group ${\Gamma_s=\Gamma\cdot\Gamma(p^s)}$ has the exact level~$p^s$ for ${1<s<e}$. (See also Remark~\ref{R.es}.)
\end{rem}

\begin{prop}\label{P.Vnontrivial}
If~$\Gamma$ contains~$-I$ and has at most two cusps,  then ${|G_2|\geqslant(p^4-p^2)/2}$ and ${V_1\neq\gen{0}}$. Under the additional assumption ${[\SL(\F_p):G_1]>2}$ we have ${\dim(V_1)\geqslant2}$.
\end{prop}
\begin{prf}
Let ${\mu_1=\mu(\Gamma_1)}$ be the index of~$G_1$ in ${\SL(\F_p)}$. Then ${|G_1|=(p^3-p)/\mu_1}$.
Since  ${\nu_\infty(\Gamma)\leqslant 2}$, Corollary~\ref{C.vinfty} and equation~\eqref{E.cardM} imply ${|G_2|\geqslant(p^4-p^2)/2}$. 
Hence 
\begin{equation*}
 |V_1|=\frac{|G_2|}{|G_1|}\geqslant\frac{(p^4-p^2)/2}{(p^3-p)/\mu_1}=p\mu_1/2.
\end{equation*}
For ${p>2}$ we have ${p\mu_1/2>1}$, while for ${p=2}$ we have ${V_1\ni I}$ by Remark~\ref{R.v_1}; in both cases ${V_1\neq\gen{O}}$.
If, in addition, ${\mu_1>2}$ then ${p\mu_1/2>p}$ and ${\dim(V_1)>1}$.
\end{prf}

\medskip

We conclude this subsection with yet another relation between the spaces~$V_s$. Although it is not explicitly used in the present article, we  include it for further references.
\begin{prop}\label{P.MN-NM}
 Let ${M_1\in V_{s_1}}$ and ${M_2\in V_{s_2}}$. Then ${M_1M_2-M_2M_1}$ lies in $V_{s_1+s_2}$.
\end{prop}
\begin{prf}
{\sloppy
By surjectivity of the reduction maps there exist matrices~$N_i$ with entries in ${\z{p^{s_j+1}}}$ (where ${\{i,j\}=\{1,2\}}$) such that ${X_i=I+p^{s_i}N_i\in G_{s_1+s_2+1}}$ projects to ${I+p^{s_i}M_i\in G_{s_i+1}}$, which means that ${N_i\equiv M_i\pmod p}$. Then over the ring $\Z/p^{s_1+s_2+1}\Z$ we have
\begin{equation*}
 X_1X_2-X_2X_1 = p^{s_1+s_2}(N_1N_2-N_2N_1),
\end{equation*}
so that the commutator of~$X_1$ and~$X_2$ is
\begin{equation*}
X_1X_2(X_2X_1)^{-1}=\left(X_2X_1+p^{s_1+s_2}(N_1N_2-N_2N_1)\right)(X_2X_1)^{-1}= I+p^{s_1+s_2}(N_1N_2-N_2N_1),
\end{equation*}
which concludes the proof.
}
\end{prf}

\medskip 

The following property will be used in Section~\ref{S.mixed}.
\begin{prop}
\label{P.pe'}
Let~$\Gamma$ be a congruence subgroup of exact level~$p^e$ and let~$\Gamma'$ be a congruence subgroup of exact level~$p^{e'}$ with ${\Gamma<\Gamma'}$. Then the index ${[\Gamma':\Gamma]}$ divides ${p^{3e-2}(p+1)(p-1)}$ and is divisible by ${p^{e-e'}}$.
\end{prop}

	\begin{prf}
	The first statement is obvious because 
	 ${|\SL(\z{p^{e}})|=p^{3e-2}(p+1)(p-1)}$.
	Further, 
	$$
\Gamma=\Gamma_e\leqslant\Gamma_{e'}\leqslant\Gamma'_{e'}=\Gamma',
$$ 
and by Corollary~\ref{C.sla}, for every  ${s<e}$ we have ${V_s\subsetneq\sla(\F_p)}$, which implies that ${[\Gamma_{s+1}: \Gamma_s]}$ is divisible by~$p$. Then~$p^{e-e'}$ divides  ${[\Gamma_{e'}:\Gamma_e]}$ and, a fortiori,  ${[\Gamma':\Gamma]}$.
	\end{prf}

\medskip

We are now ready to begin our inspection on groups of prime power level. 
We shall start with the groups such that ${p\neq2}$ divides the order of~$G_1$, then turn to those such that ${p\neq2}$ does not divide the order of~$G_1$, and finally consider the case ${p=2}$.


\subsection{The ``unipotent'' case}
Throughout this and the following subsection we shall assume ${p\neq2}$.
In this subsection we consider groups~$\Gamma$ such that~$p$ divides the order of~$G_1$. (One may call such~$\Gamma$ ``unipotent''.)

Assume that  ${\Gamma\ni -I}$ and ${\nu_\infty(\Gamma)\leqslant 2}$.  As follows from the results of Section \ref{S.prime}, the group~$G_1$ is either ${\SL(\F_p)}$ or  ${\ST(\F_p)}$, up to conjugation.  

We begin by studying the adjoint representations of ${\ST(\F_p)}$ and ${\SL(\F_p)}$, in order to find the subspaces of ${\sla(\F_p)}$ that are stable under their action.

Fix a generator~$g$ of the multiplicative group~$\F_p^\ast$ and consider the matrices 
$$
S=\mat{0}{1}{-1}{0}, \qquad T=\mat{1}{1}{0}{1}, \qquad X=\mat{g^{-1}}{0}{0}{g}
$$
in $\SL(\F_p)$.
The element~$T$ generates the maximal unipotent group ${\{\smat{1}{\ast}{}{1}\}}$; the elements~$T$ and~$X$, together, generate the special triangular group ${\ST(\F_p)}$; the three elements~$S$, $T$, and~$X$ generate the special linear group\footnote{Actually, already~$S$ and~$T$ generate ${\SL(\F_p)}$, but it is more convenient for us to include~$X$ in the set of generators.} ${\SL(\F_p)}$.

We fix for ${\sla(\F_p)}$ the basis
\begin{equation*}
A=\mat{1}{0}{0}{-1},\quad B=\mat{0}{1}{0}{0},\quad C=\mat{0}{0}{1}{0}.
\end{equation*}

\begin{prop}
\label{P.SLinvariant}
For ${p\neq2}$ the only proper non-zero  ${\ST(\F_p)}$-invariant subspaces of ${\sla(\F_p)}$ are $\gen{B}$ and $\gen{A,B}$.  There are no proper non-zero ${\SL(\F_p)}$-invariant subspaces of ${\sla(\F_p)}$.
\end{prop}

\begin{prf}
We consider for ${\sla(\F_p)}$ the basis
\begin{equation*}
e_1=4B,\quad e_2=2A,\quad e_3=-A-2C.
\end{equation*}
In this basis, the conjugation map ${M\mapsto T^{-1}MT}$ has the (left) matrix
\begin{equation*}
\begin{pmatrix}
1&1&0\\0&1&1\\0&0&1
\end{pmatrix}.
\end{equation*}
Hence the proper non-zero $T$-invariant subspaces of ${\sla(\F_p)}$ are ${\gen{e_1}=\gen{B}}$ and ${\gen{e_1,e_2}=\gen{A,B}}$. Since both are also $X$-invariant, they are ${\ST(\F_p)}$-invariant, and there are no other. Since none of them is $S$-invariant, there is no   non-zero proper invariant ${\SL(\F_p)}$-subspaces. 
\end{prf}

\medskip

This proposition allows us to settle the case ${G_1=\SL(\F_p)}$, for ${p>2}$. 
\begin{coro}\label{C.SL}
Let~$\Gamma$ be of level $p^e$ with ${p\neq2}$. Assume that ${\Gamma \ni -I}$, that ${\nu_\infty(\Gamma)\leqslant 2}$ and that ${G_1=\SL(\F_p)}$. Then ${\Gamma=\SL(\Z)}$.
\end{coro}
\begin{prf}
Propositions~\ref{P.Vinvariant},~\ref{P.Vnontrivial} and \ref{P.SLinvariant} imply ${V_1=\sla(\F_p)}$, and we conclude by Corollary~\ref{C.sla}.
\end{prf}

\medskip
Now we are ready to classify all ``unipotent''~$\Gamma$ of odd prime power level, such that
\begin{equation}
\label{propell}
\parbox{0.7\textwidth}{$\Gamma\ni -I$, and every congruence subgroup of~$\Gamma$ containing the elliptic elements of~$\Gamma$ has at most~$2$ cusps}
\end{equation}



\begin{prop}
\label{P.power_unip}
Let~$\Gamma$ be a congruence subgroup of exact level~$p^e$, with ${e>1}$ and ${p\neq2}$, such that~$p$ divides ${|G_1|}$. Assume that~$\Gamma$ satisfies~\eqref{propell}. Then 
we have one of the following two cases:
\begin{itemize}
\addtolength{\itemsep}{-0.5\baselineskip}
\item
$\Gamma$ is of exact level~$9$ and the curve~$X_\Gamma$ is of genus~$0$;

\item
$\Gamma$ is of exact level~$27$ and the curve~$X_\Gamma$ is of genus~$1$. 

\end{itemize} 
(In both cases~$\Gamma$ is uniquely defined up to conjugacy.)
\end{prop}

Together with Theorem~\ref{T.B_subgroup} this has the following consequence. 

\begin{coro}
\label{C.unip}
Let~$\Gamma$ be a congruence subgroup of exact level~$p^e$, with ${e>1}$ and ${p\neq2}$, such that~$p$ divides ${|G_1|}$. Then either Siegel's theorem is effective for~$X_\Gamma$ or~$X_\Gamma$ is non-Siegelian.
\end{coro}

\begin{rem}
One can give a totally explicit description for the groups~$\Gamma$ from Proposition~\ref{P.power_unip}. For instance, for the~$\Gamma$ of level~$9$, the group~$G_2$ is, up to conjugacy generated by the matrices
$$
\mat{2}{-1}{3}{-1},
\mat{1}{3}{0}{1},
\mat{4}{0}{0}{-2}
 \in \SL(\Z/9\Z). 
$$
One can exhibit a  similar set of generators for the group~$G_3$ which defines the~$\Gamma$ of level~$27$; see \cite[Proposition~7.3.5]{Il08} for the missing details. 
\end{rem}

\begin{prf}[of Proposition~\ref{P.power_unip}]
If~$\Gamma$ contains~$-I$ and has at most~$2$ cusps,  then so does~$\Gamma_1$, and by Theorem~\ref{T.p.list} we may assume that either ${G_1=\SL(\F_p)}$ or ${G_1=\ST(\F_p)}$. The former case is impossible by Corollary~\ref{C.SL}, so we have ${G_1=\ST(\F_p)}$. 

Next, let~$G'$ be the subgroup of $\ST(\F_p)$ generated by elements of order dividing~$12$, and let~$\Gamma'$ be the intersection of~$\Gamma$ with the pull-back of~$G'$ to $\SL(\Z)$. Then~$\Gamma'$ contains $-I$ and the elliptic elements of~$\Gamma$. 
On the other hand, for ${p\notin \{2,3,5,7,13\}}$ the group~$G'$ is a proper subgroup of $\ST(\F_p)$, and Proposition~\ref{P.st} implies that $\Gamma'_1$ has at least~$3$ cusps. Hence so does~$\Gamma'$. We conclude that ${p\in\{3,5,7,13\}}$.

The group~$\Gamma$ contains~$-I$, has at most two cusps, and satisfies ${\mu(\Gamma_1)=p+1}$.  Proposition~\ref{P.Vnontrivial} implies ${\dim(V_1)\geqslant2}$.
By Proposition~\ref{P.Vinvariant}, Proposition~\ref{P.SLinvariant} and Corollary~\ref{C.sla} we obtain ${V_1=\gen{A,B}}$, and Proposition~\ref{P.Vcontained} gives ${V_1=V_2=\ldots=V_{e-1}}$. Hence ${|G_e|=p^{2e}(1-p^{-1})}$. 

On the other hand,  Corollary~\ref{C.vinfty} and~\eqref{E.cardM} imply that any subgroup~$\Gamma'$ of~$\Gamma$ that contains~$-I$ and has at most~$2$ cusps must satisfy 
${|G'_e|\geqslant|\M_{p^e}|/2=p^{2e}(1-p^{-2})/2}$. 
Thus,
\begin{equation*}
[G_e:G'_e]\leqslant\frac{p^{2e}(1-p^{-1})}{p^{2e}(1-p^{-2})/2}=\frac{2}{1+p^{-1}}<2.
\end{equation*}
This implies that if~$\Gamma$ satisfies~\eqref{propell} then any proper subgroup of~$\Gamma$ of the exact level~$p^e$, containing~$-I$, cannot have at most two cusps. In particular, if~$\Gamma$ satisfies~\eqref{propell} then the congruence subgroup generated by ${\Gamma(p^e)}$, by~$-I$, and by the elliptic elements of~$\Gamma$ is~$\Gamma$ itself.

A direct verification on the levels $3^4$, $5^2$, $7^2$, and ${13^2}$ shows that there exist no groups~$\Gamma$ with this property, with at most two cusps, and such that ${G_1=\ST(\F_p)}$.
A further inspection on the levels~$3^2$ and~$3^3$ concludes our classification. See \cite[Section~4.3]{Il08} for more details. 

If~$\Gamma$ does not satisfy~\eqref{propell} then Siegel's theorem is effective for~$X_\Gamma$ by Theorem~\ref{T.B_subgroup}. If~$\Gamma$ satisfies~\eqref{propell} then either~$X_\Gamma$ has genus~$1$, and Siegel's theorem is effective for~$X_\Gamma$ by Theorem~\ref{T.classic}, or~$X_\Gamma$ has genus~$0$, and is non-Siegelian.
\end{prf}

\subsection{The ``semi-simple'' case}
As in the previous subsection, we assume ${p\neq2}$. In this subsection we consider groups~$\Gamma$ such that~$p$ does not divide the order of~$G_1$. (One may call such~$\Gamma$ ``semi-simple''.)
As we have seen in Section~\ref{S.prime}, up to conjugacy there are ten possible groups~$G_1$ for ${p\neq2}$.

We shall need a simple lemma, that will be used for ${n=3}$, but we state the general case. It is certainly well-known, but we include a proof for the sake of completeness.
\begin{lemm}
\label{L.basis}
Let~$A$ be an algebra over a field of characteristic distinct from~$2$, and let ${X_1, \ldots, X_n}$ be invertible and pairwise anti-commuting\footnote{that is, ${X_iX_j+X_jX_i=0}$ for ${i\ne j}$} elements of~$A$. Then ${X_1, \ldots, X_n}$ are linearly independent over the base field.
\end{lemm}
\begin{prf}
Let ${S=\sum_ia_iX_i}$ be a linear combination of the~$X_i$, with~$a_i$ in the base field. If ${S=0}$ then for every~$i$ we have
\begin{equation*}
 0=X_iS+SX_i=\sum_{j\neq i}a_i(X_iX_j+X_jX_i)+2a_iX_i^2=2a_iX_i^2.
\end{equation*}
Since~$X_i$ is invertible in~$A$ and~$2$ is invertible in the base field, this implies that every~$a_i$ is~$0$.
\end{prf}

\medskip

Now we have the following property, which allows us to exclude immediately seven of the ten cases referred to in the beginning of this subsection.
\begin{prop}\label{P.basis}
Let~$G$ be a subgroup of ${\SL(\F_p)}$ and let~$\bar G$ be its image in~$\PSL(\F_p)$.
If~$\bar G$ contains a subgroup~$\bar H$ isomorphic to~$\z{2}\times\z{2}$, then  ${\sla(\F_p)}$ has a basis consisting of three elements of~$G$.

If~$\bar G$ contains a subgroup isomorphic to the alternating group~$\Aa_4$, then there are no non-trivial $G$-stable subspaces of ${\sla(\F_p)}$.
\end{prop}
\begin{prf}
Let~$\bar X$ and~$\bar Y$ be generators of ${\bar H\equiv\z{2}\times\z{2}}$ and let ${\pm X}$ and ${\pm Y}$ be their pullbacks in~$G$. Since the elements~$X$, $Y$, and ${XY}$ are traceless, they belong to ${\sla(\F_p)}$. The obvious relations
${X^2=Y^2=(XY)^2=-I}$
show that~$X$, $Y$, and~$XY$ are pairwise anti-commuting as in Lemma~\ref{L.basis}. Hence they form a basis of $\sla(\F_p)$. 

In this basis, the conjugation maps by~$X$,~$Y$, and ${XY}$ have the matrices
\begin{align*}
 \gamma_X&=
 \begin{pmatrix}
  1&0&0\\
  0&-1&0\\
  0&0&-1
 \end{pmatrix},&
 \gamma_Y&=
 \begin{pmatrix}
  -1&0&0\\
  0&1&0\\
  0&0&-1
 \end{pmatrix},&
 \gamma_{XY}&=
 \begin{pmatrix}
  -1&0&0\\
  0&-1&0\\
  0&0&1
 \end{pmatrix}.
\end{align*}
This implies that the $G$-invariant subspaces of $\sla(\F_p)$ are generated by subsets of ${\{X,Y,XY\}}$.

Let now~$\bar G$ contain a subgroup isomorphic to~$\Aa_4$; in turn, this will contain a subgroup~$\bar H$ isomorphic to the Klein group ${T\cong\z{2}\times\z{2}}$, and an element~$\bar R$ that cycles, by conjugation, the non-trivial elements of~$\bar H$. 
Taking a basis~$X$, $Y$, and ${XY}$ of ${\sla(\F_p)}$ as above, the pullback~$R$ of~$\bar R$ in~$G$ cycles the spaces $\gen{X}$, $\gen{Y}$, and $\gen{XY}$. Thus the only $G$-invariant subspaces of ${\sla(\F_p)}$ are trivial.
\end{prf}

\begin{prop}
\label{P.power_semis}
Let~$\Gamma$ be a congruence subgroup of exact level~$p^e$, with ${e>1}$ and ${p\neq2}$, such that~$p$ does not divide ${|G_1|}$. Assume that~$\Gamma$ contains~$-I$ and has at most two cusps. Then 
we have one of the following cases:
\begin{itemize}
\addtolength{\itemsep}{-0.5\baselineskip}
\item
${p^e=9}$, ${\bar G_1=\Cc_2}$ and the curve~$X_\Gamma$ is of genus~$0$;
\item
${p^e=9}$, ${\bar G_1=\Dd_2}$ and the curve~$X_\Gamma$ is of genus~$0$;
\item
${p^e=25}$, ${\bar G_1=\Dd_3}$ and the curve~$X_\Gamma$ is of genus~$2$.
\end{itemize}
(In all three cases the group~$\Gamma$ is uniquely defined up to conjugacy.)
\end{prop}

Together with Theorem~\ref{T.B_3cusps} this has the following consequence. 

\begin{coro}
\label{C.semis}
Let~$\Gamma$ be a congruence subgroup of exact level~$p^e$, with ${e>1}$ and ${p\neq2}$, and such that~$p$ does not divide ${|G_1|}$. Then either Siegel's theorem is effective for~$X_\Gamma$, or~$X_\Gamma$ is non-Siegelian, or ${p^e=25}$. \qed
\end{coro}
 
\begin{rem}
Again, one can give a more explicit description for the three groups~$\Gamma$ above. For instance, for the ``wicked''~$\Gamma$ of level~$25$, the group~$G_2$ is, up to conjugacy, generated by the matrices
\begin{equation*}
\mat{0}{1}{-1}{0},\mat{0}{7}{7}{-1},\mat{1}{5}{-5}{1},\mat{11}{-5}{0}{-9} \in \SL(\Z/25\Z).
\end{equation*}
One has a similar description for the two groups of level~$9$; see \cite[Proposition~7.4.4]{Il08}. 
\end{rem}

\begin{prf}[of Proposition~\ref{P.power_semis}]
If~$\Gamma$ contains~$-I$ and has at most~$2$ cusps, then  then so does~$\Gamma_1$. Theorem~\ref{T.p.list} now implies that, up to conjugacy,~$G_1$ is one of the ten groups with ${p\neq2}$ in Table~\ref{Tab.p8}.
We also have ${\mu(G_1)>2}$, which implies ${\dim(V_1)\geqslant2}$ by Proposition~\ref{P.Vnontrivial}. 
Now the seven groups corresponding to the final four lines of Table~\ref{Tab.p8} can be excluded 
using Proposition~\ref{P.Vinvariant}, Proposition~\ref{P.basis} and Corollary~\ref{C.sla}.

We are left with the cases when either ${p=3}$ and ${\bar{G}_1\cong\Cc_2}$, or ${p=3}$ and ${\bar{G}_1\cong\Dd_2}$, or ${p=5}$ and ${\bar{G}_1\cong\Dd_3}$.
A direct verification on the levels~$3^3$ and~$5^3$ shows that there exist no groups~$\Gamma$ of these exact levels that contain~$-I$, have at most~$2$ cusps and such and such that~$p$ does not divide~${|G_1|}$. A further inspection on the levels~$3^2$ and~$5^2$ concludes our classification. 

We conclude the proof using Theorem~\ref{T.B_3cusps}.
\end{prf}

\subsection{The case ${p=2}$}
In this subsection we assume ${p=2}$, that is, $G_1$ is a subgroup of ${\SL(\F_2)\cong\Ss_3}$. The following propositions are proved by inspection on the levels~$2^s$ for ${s\leqslant6}$. For the details see \cite[Section~7.5]{Il08}.

\begin{prop}\label{P.2unipotent}
Let~$\Gamma$ be a congruence subgroup of exact order~$2^e$ with ${e>1}$ and with ${G_1\cong \Cc_2}$. Assume that~$\Gamma$ contains~$-I$ and satisfies~\eqref{propell}. Then ${e\leqslant5}$ and~$\Gamma$ is uniquely determined by~$e$ up to conjugacy. For ${2\leqslant e\leqslant4}$ the curve~$X_\Gamma$ is non-Siegelian, while for ${e=5}$ the curve $X_\Gamma$ is Siegelian and has genus~$1$. \qed
\end{prop}
\begin{prop}\label{P.2ST}
Let~$\Gamma$ be a congruence subgroup of exact level~$2^e$ with ${e>1}$ and with ${G_1\cong \Cc_3}$. Assume that~$\Gamma$ contains~$-I$ and satisfies~\eqref{propell}. Then ${e\leqslant4}$ and~$\Gamma$ is uniquely determined by~$e$ up to conjugacy. Moreover~$X_\Gamma$ is non-Siegelian. \qed
\end{prop}
\begin{prop}\label{P.2SL}
Let~$\Gamma$ be a congruence subgroup of exact level~$2^e$, for some ${e>1}$, containing~$-I$, having at most~$2$ cusps and with ${G_1= \SL(\F_2)}$. Then ${e\leqslant4}$ and~$\Gamma$ belongs to one of eight distinct conjugacy classes. For each of them~$X_\Gamma$ is non-Siegelian. \qed
\end{prop}

Together with Theorems~\ref{T.classic} and~\ref{T.B_subgroup}, the above results have the following consequence.
\begin{coro}
\label{C.p=2}
Let~$\Gamma$ be a congruence subgroup of exact level~$2^e$ with $e>1$. Then either Siegel's theorem is effective for~$X_\Gamma$ or~$X_\Gamma$ is non-Siegelian. \qed
\end{coro}

Theorem~\ref{T.power} is a combination of Theorem~\ref{T.prime} and Corollaries~\ref{C.unip},~\ref{C.semis} and~\ref{C.p=2}. 

The complete least of non-Siegelian curves of exact level~$p^e$ with ${e>1}$ is given in Table~\ref{Tab.power}.

\begin{table}[hbt!]
 \caption{Non-Siegelian modular curves~$X_\Gamma$ of exact prime power level~$p^e$ with ${e>1}$}
  \label{Tab.power}
  \begin{center}
  \begin{tabular}{cccccclcl}
    $p^e$ & & $\Gamma_1$ & $\mu$ & $\nu_\infty$ & $\nu_2$ & $\nu_3$ & $\g$ & remark \\
    \hline
     4 & \multirow{3}{*}{$\left.\rule{0ex}{4.2ex}\right\}$} & \multirow{3}{*}{$G_1\cong\Cc_2$} &
            6 & 2 &  2 &  0 & 0 & \\
     8 &&& 12 & 2 &  4 &  0 & 0 & \\
    16 &&& 24 & 2 &  8 &  0 & 0 & \\
     4 & \multirow{3}{*}{$\left.\rule{0ex}{4.2ex}\right\}$} & \multirow{3}{*}{$G_1\cong\Cc_3$} & 
            8 & 2 &  0 &  2 & 0 & \\
     8 &&& 16 & 2 &  0 &  4 & 0 & \\
    16 &&& 32 & 2 &  0 &  8 & 0 & \\
     4 & \multirow{4}{*}{$\left.\rule{0ex}{5.6ex}\right\}$} & \multirow{4}{*}{$\Gamma_1=\Gamma(1)$} & 
            4 & 1 &  2 &  1 & 0 & \\
     8 &&& 16 & 2 &  4 &  1 & 0 & \\
     8 &&&  8 & 1 &  2 &  2 & 0 & 2 groups\\
    16 &&& 16 & 1 &  2 &  4 & 0 & 4 groups\\
     9 && $\Gamma_1=\Gamma_0(3)$ &
           12 & 2 &  0 &  3 & 0 & \\
     9 && $\bar{G}_1\cong\Cc_2$ & 
           18 & 2 &  6 &  0 & 0 & \\
     9 && $\bar{G}_1\cong\Dd_2$ & 
            9 & 1 &  5 &  0 & 0 &
  \end{tabular}
  \end{center}
\end{table}

\section{The Mixed Levels}
\label{S.mixed}

In this section we study groups of mixed level. Our goal is to prove Theorem~\ref{T.mixed}. 
Let~$\Gamma$ be a congruence subgroup of exact level~$n$, and let the factorization of~$n$ be 
$$
n=\prod_{i\in I}q_i=\prod_{i\in I}p_i^{e_i}, 
$$
where the~$p_i$ are distinct primes and ${e_i>0}$ for every ${i\in I}$.
For every positive integer~$d$ we denote by~$\Gamma_d$ the composite group ${\Gamma\cdot\Gamma(d)}$, of level dividing~$d$, and by ${\G_d<\SL(\z{d})}$ its projection modulo~$d$. 
The group ${\SL(\z{n})}$ is isomorphic to the direct product ${\prod_{i\in I}\SL(\z{q_i})}$; this allows us to consider ${\G=\G_n}$ as a subgroup of the direct product ${\prod_{i\in I}\G_{q_i}}$.

\begin{rem}
\begin{enumerate}
\addtolength{\itemsep}{-0.5\baselineskip}
\item
If ${n=p^e}$ then~$G_s$ and~$\Gamma_s$ of Section~\ref{S.power} become~$G_{p^s}$ and~$\Gamma_{p^s}$ in this section. 
\item
If Siegel's theorem is effective for~${\Gamma_d}$ then it is effective also for~$\Gamma$.
\item
Notice that for ${d\neq n}$ the group~$\Gamma_d$ is not obliged to have the exact level~$d$, even if~$d$ divides~$n$; as in Remark~\ref{R.es}, {\sl a priori} we only know that~$\Gamma_d$ is of level dividing~$d$.
\end{enumerate}
\end{rem}

\subsection{Proof of Theorem~\ref{T.mixed}}
We begin with the following useful observation. Let ${\{S_i\}_{i\in I}}$ be a finite family of finite groups~$S_i$ and  let ${S=\prod_{i\in I} S_i}$ be their direct product. 
For a subset ${J\subset I}$ we view ${S_J=\prod_{i\in J} S_i}$ as a subgroup of~$S$, and we denote by
${\pi_J\colon S\to S_J}$ the natural projection.
\begin{prop}
Let~$T$ be a subgroup of~$S$, and let~$T_J$ and~$U_J$ be the subgroups of~$S_J$ defined by 
${T_J=\pi_J(T)}$ and  ${U_J=T\cap S_J}$. 
Then~$U_J$ is a normal subgroup of~$T_J$. Let also~$r_i$ be the index of $U_{\{i\}}$ in $T_{\{i\}}$. Then~$r_j$ divides ${\prod_{i\neq j} r_i}$ for every ${j\in I}$.
\end{prop}
\begin{prf}
Let ${I=J\cup K}$ be a partition of~$I$. 
The group ${U_J=\Ker(\pi_K|_T)}$ is normal in~$T$; then ${U_J=\pi_J(U_J)}$ is a normal subgroup of ${T_J=\pi_J(T)}$. 
The composite map ${T\twoheadrightarrow T_J\twoheadrightarrow T_J/U_J}$ has kernel ${U_J\times U_K}$ and induces an isomorphism ${T/(U_J\times U_K)\cong T_J/U_J}$, which proves ${T_J/U_J\cong T_K/U_K}$.

Now note that 
\begin{equation*}
\prod_{i\in K} U_{\{i\}}<U_K<T_K<\prod_{i\in K} T_{\{i\}}.
\end{equation*}
{\sloppy
This implies that ${|T_K/U_K|}$ divides ${\left|\prod_{i\in K}T_{\{i\}}\big/\prod_{i\in K}U_{\{i\}}\right|=\prod_{i\in K}r_i}.\ $  Taking ${J=\{j\}}$, we obtain ${T_{\{j\}}/U_{\{j\}}\cong T_K/U_K}$, whence the result.
}
\end{prf}

\medskip

Applying the above proposition to the group ${\G<\prod_{i\in I}\G_{q_i}}$ we obtain the following.

\begin{coro}\label{C.r_i}
Let~$\Gamma$ be a congruence subgroup of exact level ${n=\prod_{i\in I}q_i}$. Then for every ${i\in I}$ the congruence subgroup ${(\Gamma\cap\Gamma(n/q_i))\cdot\Gamma(q_i)}$ of exact level~$q_i$ projects modulo~$q_i$ onto a normal subgroup~$\Hh_{q_i}$ of~$\G_{q_i}$ of index~$r_i$, and~$r_j$ divides ${\prod_{i\neq j}r_i}$ for every ${j\in I}$.\qed
\end{coro}

The next statement  is certainly well-known, but we include a proof for the sake of completeness.
\begin{prop}
\label{P.p}
Let~$p$ be a prime and let~$H_s$ be a normal subgroup of ${\SL(\z{p^s})}$ for some  ${s>0}$. If ${H_s\neq\SL(\z{p^s})}$ then~$p$ divides the index of~$H_s$.
\end{prop}
\begin{prf}
When ${s=1}$, the cases ${p\leqslant 3}$ are verified by inspection, and for ${p\geqslant 5}$ any proper normal subgroup of $\SL(\Z/p\Z)$ is contained in  $\{\pm I\}$. For ${s>1}$, the projection of~$H_s$ modulo~$p^{s-1}$ is a normal subgroup~$H_{s-1}$ of ${\SL(\z{p^{s-1}})}$, and
${[\SL(\z{p^s}):H_s]=p^a[\SL(\z{p^{s-1}}):H_{s-1}]}$
for some ${a\geqslant0}$. We conclude by induction.
\end{prf}

\medskip

This has the following consequence.
\begin{prop}
\label{P.p_max}
Let~$\Gamma$ be congruence subgroup of exact level~$n$, let ${p>3}$ be the largest prime divisor of~$n$ and let~$q$ be the exact power of~$p$ dividing~$n$. Then ${\G_q\neq\SL(\Z/q\Z)}$.
\end{prop}
\begin{prf}
Let $\G_q$ and $\Hh_q$ be as in Corollary~\ref{C.r_i}.
Since~$p$ does not divide ${\bigl|\SL(\z{p_i^{e_i}})\bigr|}$ for any prime ${p_i<p}$, it cannot divide ${[\G_q:\Hh_q]}$ by Corollary~\ref{C.r_i}. Proposition~\ref{P.p} now implies that if ${\G_q=\SL(\Z/q\Z)}$ then ${\Hh_q=\G_q}$, but in this case~$p$ would not divide~$n$.
\end{prf}
\begin{coro}
\label{C.p_max}
Let~$\Gamma$ and~$p$ be as in Proposition~\ref{P.p_max}. Assume that ${p>13}$. Then Siegel's theorem is effective for~$X_\Gamma$.
\end{coro}
\begin{prf}
As above, let~$q$ be the exact power of~$p$ dividing~$n$, and consider the congruence subgroup ${\Gamma_q=\Gamma\cdot\Gamma(q)}$ of level dividing~$q$. Since ${p>13}$, the results of the previous sections imply that either Siegel's theorem is effective for~$X_{\Gamma_q}$ or ${\G_q=\SL(\Z/q\Z)}$, which contradicts Proposition~\ref{P.p_max}.
\end{prf}

\medskip

\begin{prf}[of Theorem~\ref{T.mixed}]
Let~$\Gamma$ be a subgroup of exact level ${n=\prod p^{e_p}}$. If the set of prime divisors of~$n$ is not contained in ${\{2,3,5,7,11,13\}}$ then we conclude by Corollary~\ref{C.p_max}. 
Assume now that~$n$ factors in the primes $2, 3, 5, 7, 11, 13$.

For every prime~$p$ let ${\Gamma_{(p)}=(\Gamma\cap\Gamma(n/p^{e_p}))\cdot\Gamma(p^{e_p})}$ and ${\Gamma_{(p)}'=\Gamma_{p^{e_p}}=\Gamma\cdot\Gamma(p^{e_p})}$ be the congruence subgroups of Corollary~\ref{C.r_i}, of exact levels respectively~$p^{e_p}$ and~$p^{e_p'}$, and with ${\Gamma_{(p)}<\Gamma_{(p)}'}$. Put also ${r_p=[\Gamma_{(p)}':\Gamma_{(p)}]}$.

If Siegel's theorem is effective for~$X_{\Gamma'_{(p)}}$ then it is effective for~$X_\Gamma$, too. Otherwise, by the results of the previous sections, we have ${e_2'\leqslant4}$, ${e_3', e_5'\leqslant2}$, and ${e_7', e_{11}', e_{13}'\leqslant1}$.

We are now going to find a bound for~$e_2$. By Corollary~\ref{C.r_i} we have that~$r_2$ divides ${\prod_{p\neq2}r_p}$. By Proposition~\ref{P.pe'} we obtain that $2^{e_2-e_2'}$ divides~$r_2$ and that ${\prod_{p\neq2}r_p}$ divides ${\prod_{p\neq2}p^{3e_p-2}(p+1)(p-1)}$. Thus $2^{e_2-e_2'}$ divides 
${\prod_{p\leqslant 13}(p+1)(p-1)}$
which implies ${e_2-e_2'\leqslant 16}$. Since ${e_2'\leqslant4}$, we obtain ${e_2\leqslant20}$.

In exactly the same way we bound the other exponents~$e_p$, completing thereby the proof.
\end{prf}

\end{document}